\documentclass{amsart}
\usepackage[utf8]{inputenc}
\usepackage{verbatim}

\usepackage{amssymb}
\usepackage{amscd}

\usepackage{soul,xcolor}
\setstcolor{red}

\usepackage{latexsym,amssymb,amsmath}
\usepackage{epsfig}
\input xy
\xyoption{all}
\usepackage[english]{babel}

\usepackage{tikz}

\usepackage{amsmath}
\usepackage{bbm} 
\newtheorem{theorem}{Theorem}[section]
\newtheorem{lemma}[theorem]{Lemma}
\newtheorem{prop}[theorem]{Proposition}
\newtheorem{coro}[theorem]{Corollary}

\def\PP{\mathbb{P}}
\def\CC{\mathbb{C}}
\def\ZZ{\mathbb{Z}}\def\QQ{\mathbb{Q}}

\def\cC{{\mathcal C}}

\def\cQ{{\mathcal Q}}

\def\cO{{\mathcal O}}
\def\cS{{\mathcal S}}
\def\cU{{\mathcal U}}
\def\cV{{\mathcal V}}

\def\ra{\rightarrow}\def\lra{\longrightarrow}

\def\fsl{{\mathfrak{sl}}}

\def\fe{{\mathfrak e}}

\newcommand{\of}{\mathcal{O}}

\usepackage{xcolor}

\title{A four-dimensional cousin of the Segre cubic}
\author{Laurent Manivel}

\address{Institut de Math\'ematiques de Toulouse ; UMR 5219, Universit\'e de Toulouse \& CNRS, F-31062 Toulouse Cedex 9, France}
\email{manivel@math.cnrs.fr}

\begin{document}

\begin{abstract}
This note is devoted to a special Fano fourfold defined by a four-dimensional space of skew-symmetric 
forms in five variables. This fourfold appears to be closely related with the classical Segre cubic and 
its Cremona-Richmond configuration of planes. Among other exceptional properties, it is infinitesimally 
rigid and has Picard number six. We show how to construct it by blow-up and contraction, 
starting from a configuration of five planes in a four-dimensional quadric, compatibly with the 
symmetry group $\cS_5$. From this construction we are able to describe the Chow ring explicitely. 
\end{abstract}

\maketitle

\vspace{-5mm}
\centerline{{\it Dedicated to the memory of Laurent Gruson}} 

\section{Introduction}

Fano threefolds were classified more that fourty years ago, after some fifty years of efforts. The 
classification of Fano fourfolds is still elusive and will probably remain so for a long time.
There are many ways to construct such manifolds, and a systematic study was launched a few years ago,
of those that can be constructed from vector bundles on products of Grassmannians and more general 
flag manifolds \cite{bftR}; a sample has already appeared in \cite{bfm}. In this database, there is 
a unique fourfold with maximal Picard number, equal to six: the study of this fourfold is the 
object of this note.

This study turned out to be  related with interesting questions at the intersection of algebraic 
geometry with Lie theory.
Consider two complex vector spaces $V_4$ and $V_5$, of dimension four and five respectively. 
The action of $GL(V_4)\times GL(V_5)$ on $V_4^\vee\otimes\wedge^2V_5^\vee$ is known to be prehomogenous, 
its open orbit being the complement of a degree $40$ hypersurface \cite[p.98]{sk}. It is in fact one of the most 
complicated prehomogeneous spaces, containing no less than $63$ distinct orbits \cite{ozeki79, djoko06}. An important 
literature has been devoted to this prehomogeneous space, including some in connections with 
quintic field extensions, in the spirit of Bhargava's work on higher reciprocity laws \cite{kable-yukie04a, kable-yukie04b, bhargava08}. 

The Fano fourfold $X_4$ we are interested in is defined by a generic element of the prehomogeneous space $V_4^\vee\otimes\wedge^2V_5^\vee$. It has two natural projections to $G(2,4)\simeq\QQ^4$ and to the 
six-dimensional $G(3,5)$ that we describe in some details in section 4. In particular 
we show it is a small resolution of a fourfold with ten singular points which appears to be a cousin, 
or a big brother of
the Segre cubic primal; this small resolution contracts ten planes which can be seen as a special subcollection 
of the classical Cremona-Richmond configuration. We deduce:

\medskip\noindent {\bf Theorem.} {\it 
Consider five general planes in one of the two families of projective planes in $\QQ^4$. They intersect 
pairwise in ten points. Blow-up these ten points and then the strict transforms of the five planes. Then the 
strict transforms of the exceptional divisors of the first blowup can be contracted to  a smooth Fano fourfold, 
which is precisely $X_4$. }

\medskip
Then we show that the automorphism group is $Aut(X_4)=\cS_5$, so that $$ Pic(X_4)^{\cS_5}\simeq \ZZ^2$$ is generated by the
pull-back of the hyperplane classes by the two projections. This suggests to construct the tensor that 
defines $X_4$ by {\it reverse-engineering}, starting from the representation theory of $\cS_5$; we show how
this leads to a normal form from this tensor. 
We then use the previous constructions to describe the Chow ring of $X_4$ completely, including the action of $\cS_5$. 
We also check that $X_4$, as expected, is infinitesimally rigid. 

\smallskip 
This study can be considered as a warm-up for the more mysterious case of $U_5^\vee\otimes\wedge^2V_5^\vee$,
directly related to $E_8$, which has infinitely many but well-described orbits for the action of 
$GL(U_5)\times GL(V_5)$ (see \cite{glw} for a first approach). Among other nice geometric objects, this representation 
will give rise to an interesting family of special Fano sixfolds. 

\medskip {\it Acknowledgements.} We thank Marcello Bernardara, Enrico Fatighenti and Fabio Tanturri for sharing 
our joint project on Fano fourfolds. Thanks also to Pieter Belmans and Igor Dolgachev for their comments and 
suggestions. We acknowledge support from the ANR project FanoHK, grant ANR-20-CE40-0023. 

\section{Models}

According to the classical  Borel-Weil theorem, one can interprete the  representation $V_4^\vee\otimes\wedge^2V_5^\vee$
as a space of global sections of an irreducible homogeneous vector bundle
over a homogeneous space, and this in more than one way:
$$\begin{array}{rcl} 
V_4^\vee\otimes\wedge^2V_5^\vee & = & \Gamma(G(2,V_4)\times\PP(V_5), \cU^\vee\boxtimes \cQ^\vee(1)) \\
& = & \Gamma(\PP(V_4)\times\PP(V_5^\vee),\of(1)\boxtimes \wedge^2\cV^\vee) \\
& = & \Gamma(G(2,V_5),V_4^\vee\otimes \wedge^2\cV^\vee) \\ 
& = & \Gamma(\PP(V_4)\times\PP(V_5),\of(1)\boxtimes \cQ^\vee(1)) \\
 & = & \Gamma (G(2,V_4)\times G(3,V_5), \cU^\vee\boxtimes \wedge^2\cV^\vee) \\
& = & \Gamma(\PP(V_4)\times G(3,V_5),\of(1)\boxtimes \wedge^2\cV^\vee) \\ 
& = & \Gamma(\PP(V_4)\times G(2,V_5),\of(1)\boxtimes \wedge^2\cV^\vee) \\ 
 & = & \Gamma (G(2,V_4)\times G(2,V_5), \cU^\vee\boxtimes \wedge^2\cV^\vee).
 \end{array}$$
Here $\cU$ and $\cV$ denote tautological bundles on Grassmannians (with some abuse of notations since we use the 
these symbols several times for distinct bundles on different Grassmannians).
As a consequence, consider a general element $\theta$ in  $V_4^\vee\otimes\wedge^2V_5^\vee$. Interpreting 
it as a global section of a vector bundle in these seven different ways, we obtain smooth subvarieties of codimensions
equal to the ranks of the vector bundles in question, that we respectively denote as follows 
(the notation is such that $X_d$ has dimension $d$):
$$\begin{array}{lll}
X_0\subset G(2,V_4)\times\PP(V_5), & \qquad  & X_1\subset \PP(V_4)\times\PP(V_5^\vee), \\
X_2\subset  G(2,V_5), & \qquad & X_3\subset \PP(V_4)\times\PP(V_5), \\
X_4\subset  G(2,V_4)\times G(3,V_5), &\qquad & X_6\subset \PP(V_4)\times G(3,V_5),\\
X_8\subset  \PP(V_4)\times G(2,V_5), &\qquad & X'_8\subset G(2,V_4)\times G(2,V_5).
\end{array}$$

\medskip
Another obvious thing to do is to consider $\theta$ as a general
morphism from $V_4$ to $\wedge^2V_5^\vee$. The image of $\PP(V_4)$ inside $\PP(\wedge^2V_5^\vee)$ is then a generic 
projective three-plane, that has to meet the Grassmannian $G(2,V_5^\vee)$ along a set $Y_0$ of five reduced points (the degree of the 
Grassmannian being equal to five). Correspondingly, we get a set $P_0$ of five points in $\PP(V_4)$, and a set 
 $\Pi_0$ of five planes in $\PP(V_5)$, all in general position. Concretely, if we choose a basis $e_1,\ldots, e_4$ of $V_4$, with dual basis  $e_1^\vee,\ldots, e_4^\vee$ of $V_4^\vee$ and 
 decompose $\theta$ accordingly as 
$$\theta=e_1^\vee\otimes\theta_1+e_2^\vee\otimes\theta_2+e_3^\vee\otimes\theta_3+e_4^\vee\otimes\theta_4$$
then the contraction $\theta(v)=v_1\theta_1+v_2\theta_2+v_3\theta_3+v_4\theta_4$ 
 has rank two  when $[v]$  belongs to $P_0$; that is, it decomposes as $f_1^\vee\wedge f_2^\vee$ for two linear
 forms  $f_1^\vee, f_2^\vee$ whose kernels intersect along the corresponding plane in $\PP(V_5)$. 
 We will denote the five two-forms of rank two (defined up to scalars) obtained by contracting $\theta$ as $\omega_1,\ldots , \omega_5$.
 It would be natural then to  impose the normalization 
 $\omega_1+\cdots +\omega_5=0$, and decompose $\theta$ as 
 $$\theta=u_1^\vee\otimes\omega_1+u_2^\vee\otimes\omega_2+u_3^\vee\otimes\omega_3+u_4^\vee\otimes\omega_4
 +u_5^\vee\otimes\omega_5$$
 for some linear forms $u_1^\vee,\ldots , u_5^\vee$ such that $u_1^\vee +\cdots +u_5^\vee=0$.
 
 \medskip\noindent {\it Notations}.
 
 $P_0=\{p_1,\ldots , p_5\}$ is a set of five points in $\PP(V_4)$, in natural bijection with the set 
 $\{\omega_1,\ldots , \omega_5\}$, of five decomposable two-forms in $\wedge^2V_5^\vee$, that 
 define five points in $G(2,V_5^\vee)\simeq G(3,V_5)$, hence five planes $P_1,\ldots ,P_5$ in $\PP(V_5)$. 
 They also define five planes $\pi_1,\ldots, \pi_5$ in $G(2,V_4)$, where $\pi_k$ is the set of planes in $V_4$ that contain $p_k$. 

 $L_0$ is the set of pairs of points in $P_0$. According to the previous identifications, it is in natural 
 bijection with a set of ten lines in $\PP(V_4)$, a set of ten points in $\PP(V_5)$, and a set of ten points 
 in $G(2,V_4)$. 
 
 \section{Small dimensions}
 
Most results in this section are classical. Our purpose is mainly to set up the scene for the main character, 
which will make its entry in the next section. 

\begin{prop}
$X_0$ consists in $10$ points of $G(2,V_4)\times\PP(V_5)$, in natural bijection with $L_0$. 
\end{prop}


\proof By definition, a point  $(A_2,B_1)$ belongs to $X_0$ if and only if we can decompose $\theta$ 
in such a way that $A_2$ is cut out by the linear forms $e_3^\vee, e_4^\vee$ and the skew-symmetric forms 
$\theta_1,\theta_2$ have the same kernel $B_1\subset V_5$. Otherwise said, $\theta_1$ and $\theta_2$ belong 
to $\wedge^2B_1^\perp$. Since in the latter space, decomposable tensors are parametrized by a quadric, 
we can make a change of basis in $A_2^\perp$ and suppose that  $\theta_1$ and $\theta_2$ are indeed 
decomposable. Concretely, this means that we can write $\theta$ in the form
$$\theta=e_1^\vee\otimes f_1^\vee\wedge f_2^\vee+e_2^\vee\otimes f_3^\vee\wedge f_4^\vee+e_3^\vee\otimes\theta_3+e_4^\vee\otimes\theta_4.$$
Then $[e_1]$ belongs to $P_0$, the associated plane in $\PP(V_5)$ being $P_1=\langle f_1,f_2\rangle^\perp$, and also 
$[e_2]$ belongs to $P_0$, the associated plane being $P_2=\langle f_3,f_4\rangle^\perp$. In particular $A_2=\langle e_1,e_2\rangle$ and 
$B_1=P_1\cap P_2$, as claimed. \qed


\begin{prop} 
$X_1$ is the union of five disjoint lines, in  natural bijection with $P_0$. 
\end{prop}

\proof By definition, a point in $X_1$ is a pair $(A_1,B_4)$ such that $\theta(v)$ vanishes on $B_4$ when 
$v$ generates $A_1$. But then $\theta(v)$ must have rank two, of the form $f_1^\vee\wedge f_2^\vee$. In 
particular $A_1$ must correspond to one of the five points of $P_0$, and the hyperplane $B_4$ can move in the
pencil $\langle f_1^\vee, f_2^\vee\rangle$. \qed 


\begin{prop} 
$X_2\subset G(2,V_5)$ is a del Pezzo surface of degree five. 
\end{prop}

\proof Obvious. \qed 

\medskip Recall that the del Pezzo surface of degree five contains $10$ lines. Since the embedding in 
$G(2,V_5)$ is anticanonical, this means in our setting that there exists ten flags $A_1\subset A_3\subset V_5$
such that $\theta(v,w)=0$ for any $v\in A_1, w\in A_3$. It is easy to see that these ten flags are in natural
bijection with $L_0$, the ten points $[A_1]$ in $\PP(V_5)$ being exactly the intersections of the planes $P_1,
\ldots , P_5$. 

\begin{prop} 
The projection of $X_3$ to $\PP(V_4)$ is the blow-up of the five points of $P_0$.
The projection to $\PP(V_5)$ is a small resolution of a Segre cubic primal $C_3$, ten lines being 
contracted to the ten singular points of $C_3$ defined by $L_0$.
\end{prop}

\proof (Well known.) By definition, $X_3$ parametrizes the pairs $(A_1=[v],B_1)$ such that $B_1$ is contained in the kernel 
of $\theta(v)$. Generically this two-form has rank four and the kernel is one-dimensional, which implies that 
$X_3$ projects birationally to $\PP(V_4)$. The projection has non trivial fibers when the rank of $\theta(v)$
drops, that is, over one of the five points in $P_0$. Then the kernel has dimension three and the fiber is 
a projective plane, as it has to be. 

No we turn to the projection to $\PP(V_5)$. By definition, the fibers are linear subspaces defined by the 
image of the morphism $Q(-1)\ra V_4^\vee\otimes\of_{\PP(V_5)}$ induced by $\theta$. In particular the 
fibers are non trivial over the corresponding determinantal locus $C_3$, which is a cubic threefold since 
$\det(Q(-1))=\of_{\PP(V_5)}(-3)$. This threefold becomes singular exactly when the rank drops to two. 
If $w\in V_5$ generates $B_1$, this means that the morphism from $V_4$ to $V_5^\vee$ sending $e_i$ to 
$\theta_i(w,\bullet)$ has rank two. So we may suppose after a change of basis that $\theta_1(w,\bullet)
=\theta_2(w,\bullet)=0$. In other words, $\theta_1$ and $\theta_2$ have the same kernel $B_1$, and after 
another change of basis if necessary we have already seen that we can suppose they are decomposable. 
So they define two points in $P_0$, in such a way that $B_1$ is the point obtained as the intersection of 
the corresponding planes in $\PP(V_5)$, while the line contracted to this point is the span of the 
corresponding points in $\PP(V_4)$. \qed 

As a result, $C_3$ is a cubic threefold with $10$ nodes.  (In fact $C_3$ is the image of the rational map 
from $\PP(V_4)$ to $\PP(V_5)$ sending $[v]$ to the kernel of the two-form $\theta(v)$, and essentially 
by definition this is  a Segre cubic primal  \cite{dolgachev}.

\medskip\noindent {\bf Reminder on the Segre cubic primal}. 
Recall that the Segre primal can be defined, if $x_0,\ldots , x_5$ are homogeneous coordinates on $\PP^5$, by the two equations
$$x_0+\cdots +x_5=0, \qquad x_0^3+\cdots +x_5^3=0.$$
This presentation exhibits an $\mathcal{S}_6$ symmetry, and it is known that $Aut(C_3)=\mathcal{S}_6$. 
Classically, the Segre primal contains $15$ planes. (See \cite[Chapter 9]{CAG} for much more information.)

\smallskip
The Segre cubic primal 
admits a classical modular interpretation,  according to which 
$C_3\simeq (\PP^1)^6/\hspace*{-1mm}/ SL_2$. Moreover $\bar{M}_{0,6}$ is a resolution of its
singularities (that just blows-up the singular points) and according to Kapranov it can be constructed by blowing-up five 
general points in $\PP^3$, plus the strict transforms of the ten lines that join
them \cite{kap}. 
(Note also that $\bar{M}_{0,6}$ compactifies the moduli space of genus $2$ curves.) 


\smallskip Note also that $C_3$ is known to be $G$-birationally
rigid, and even $G$-birationally superrigid, when $\mathcal{A}_5\subset G\subset \mathcal{S}_6$ \cite{avilov}. 

\smallskip
Blowing-up the ten singular points in $C_3$ we get ten exceptional divisors 
isomorphic to $\PP^1\times\PP^1$, each of which is contracted to $\PP^1$ in $X_3$. According to \cite{fink} any of the rulings of  $\PP^1\times\PP^1$ can in fact be contracted, yielding $2^{10}=1024$ small resolutions of the singularities of $C_3$, falling into $13$ orbits of  $\mathcal{S}_6$, including $6$ for which the resolution is projective. 
Homological Projective Duality for the Segre cubic is discussed in \cite{bb22}.

\medskip\noindent {\bf On the planes in the Segre cubic}. In coordinates, the $15$ planes on the Segre cubic are given by three equations 
$$x_a+x_b=x_c+x_d=x_e+x_f=0,$$
for $(abcdef)$ a permutation of $(123456)$; we denote such a plane by $(ab|cd|ef)$. Together with the $15$ points in the hyperplane $x_0+\cdots +x_5=0$ with four coordinates equal to zero, they form a $(15_3, 15_3)$ configuration classically known as the Cremona-Richmond configuration: each plane contains three of the $15$ points and each of those points belongs 
to three planes of the configuration. But beware that two planes may meet along a single point, or a projective line; the second possibility occurs when their symbols have a common pair.

\begin{prop} \label{subconfigurations}
There are exactly $6$ collections of five planes among the fifteen planes in $C_3$, meeting 
pairwise along single points. These collections are exchanged transitively by the action of $\cS_6$. 
Each one has for stabilizer a copy of $\cS_5$, embedded in $\cS_6$ in a non standard way.
\end{prop}

To understand the last sentence, recall that $\cS_6$ has the exceptional property that its outer  automorphism
group is non trivial: there exists a unique outer automorphism, and a non standard embedding of $\cS_5$ in $\cS_6$
is the composition of a standard embedding by such an outer automorphism. Note that this outer automorphism of 
$\cS_6$ exchanges the two conjugacy classes consisting of transpositions on one hand, and products of three disjoint
transpositions on the other hand; the former corresponds to points, the latter to planes in the Cremona-Richmond
configuration, which is for this reason self-dual.

\proof 
Suppose given a collection of five planes, any two of which meet at a single point. This means that each plane
is represented by three pairs, none of which being shared with another plane. So we have a total amount of $15$ 
distinct pairs; necessarily, all the $15$ pairs of integers from $1$ to $6$ must appear exactly once. 

Up to permutation, we can assume that one of our planes is $(12|34|56)$. 
Then the plane containing $(13)$ is either $(13|25|46)$ or $(13|26|45)$ and up to permuting $5$ and $6$ we can suppose it is the first one. Then the other planes are determined. For example, for the one containing $(14)$, we must split $(2356)$ into 
two pairs, and since $(25)$ and $(56)$ have already been used the only possibility is $(14|26|35)$. This also shows that 
we have three choices for the plane containing $(12)$, then two choices for the plane containing $(13)$, and then 
no more choices; this means there are exactly six possibilities. Explicitely, they are the following:

$$\begin{array}{cccccc}
(12|34|56) & (12|34|56) & (12|35|46) & (12|35|46) & (12|36|45) & (12|36|45)\\
(13|25|46) & (13|26|45) & (13|24|56) & (13|26|45) & (13|25|46) & (13|24|56)\\
(14|26|35) & (14|25|36) & (14|25|36) & (14|23|56) & (14|23|56) & (14|26|35)\\
(15|24|36) & (15|23|46) & (15|26|34) & (15|24|36) & (15|26|34) & (15|23|46)\\
(16|23|45) & (16|24|35) & (16|23|45) & (16|25|34) & (16|24|35) & (16|25|34)
\end{array}
$$

\medskip
Let us denote these six configurations by $ABCDEF$. The action of $\cS_6$ on them induces a morphism 
$\cS_6\ra\cS_6$, and a direct examination shows that it sends the 
transposition $(12)$ to the permutation $(AB)(CD)(EF)$. So it has to correspond to the outer automorphism 
of $\cS_6$, and our final claim follows.
\qed 

\medskip\noindent {\it Question.} Is there an interpretation in terms of the root system $E_7$? In fact 
the Lie algebra $\fe_7$ admits a $\ZZ_3$-grading of the form
$$\fe_7 = \fsl_3\times\fsl_6\oplus (\CC^3\otimes\wedge^2\CC^6)\oplus (\CC^3\otimes\wedge^2\CC^6)^\vee,$$
and roots of $\fe_7$ defined by weights of $\CC^3\otimes\wedge^2\CC^6$ can be interpreted as triples 
of pairs \cite{man-config}. Note that roots of $\fe_7$ are classically connected with the $28$ bitangents of 
a plane quartic.

\section{The Fano fourfold} 

Recall that our main character $X_4\subset  G(2,V_4)\times G(3,V_5)$ is defined by $\theta$ a general 
element in $V_4^\vee\otimes\wedge^2V_5^\vee$, considered as a general section of the vector bundle 
$\cU^\vee\boxtimes \wedge^2\cV^\vee$. Here $\cU$ denotes the tautological rank two bundle on $G(2,V_4)$, 
while $\cV$ denotes the tautological rank three bundle on $G(3,V_5)$.

In this section we describe the geometry of $X_4$ by blowups and contractions. 

\subsection{The main invariants} 

We start by computing the main numerical invariants of $X_4$, including its Hodge numbers.

\begin{prop}\label{invariants}
$X_4$ is a rational Fano fourfold of index one. 

Its cohomology is pure, with  $h^{1,1}=6$ and  $h^{2,2}=17$.

Moreover $h^0(-K_{X_4})=40$ and $K_{X_4}^4=172$. 
\end{prop}

\proof The first assertion is an immediate consequence of the adjunction formula. 

The Hodge numbers and invariants can computed using exact sequences, along the lines explained in \cite{bfm}.
(They could also be deduced from the geometric descriptions that will follow.)
Since $172=4\times 43$ is not divisible by 
any fourth power, the index must be one. \qed 

\medskip
Note that $h^0(-K_{X_4})=40< \dim (\wedge^2V_4\otimes\wedge^3V_5)=60$, which means that $X_4$ is linearly degenerate inside 
$G:=G(2,V_4)\times G(3,V_5)$. This can be checked by  considering the twisted Koszul complex
$$0\lra \wedge^6 E^\vee(1,1) \lra \cdots \lra E^\vee(1,1)\lra \of_{G}(1,1)\lra\of_{X_4}(1)\lra 0.$$
Indeed $H^0(E^\vee(1,1))\simeq V_4\otimes V_5$ has dimension $20$, while it can be checked that 
$H^0(\wedge^k E^\vee(1,1))=0$ for $k>1$. 

\medskip We will describe in some details the two projections $p_1, p_2$:
$$\xymatrix{&  X_4 \ar[dl]_{p_1}\ar[dr]^{p_2} & \\
G(2,V_4) & & G(3,V_5). 
}$$
We start with the second one. 

\subsection{The second projection and the Cremona-Richmond configuration} 

We start with the projection to $G(3,V_5)$, which is very similar to the resolution of singularities 
of the Segre cuic primal. 

\begin{prop}\label{p2}
The projection of $X_4$ to $G(3,V_5)$ is a small resolution of a codimension two subvariety $C_4$ of degree $12$, contracting ten planes to ten singular points in natural bijection with $L_0$. 
\end{prop}

\proof
The fiber of $p_2: X_4\lra G(3,V_5)$ 
over a point $[V]\in G(3,V_5)$ is defined by the morphism $\theta_V: \wedge^2V\ra V_4^\vee$ induced by $\theta$. In particular
the fibers are non trivial when the rank is at most two, which happens in codimension two. We conclude that the image $C_4$ of $X_4$ is a 
determinantal fourfold. Its structure sheaf admits a resolution by the
Lascoux complex \cite{lascoux}
\begin{equation}\label{lascoux}
0\lra \mathcal{V}^\vee(-3) \lra V_4\otimes \of_{G(2,V_5)}(-2)\lra \of_{G(2,V_5)}\lra\of_{C_4}\lra 0,
\end{equation}
where $\mathcal{V}$ denotes the rank three tautological bundle. 
We deduce in particular that the class of $C_4$ in the Chow ring 
of the Grassmannian $G(3,V_5)$ is $3\sigma_{11}+2\sigma_{2}$, so that its degree is $3\times 2+2\times 3=12$. 

The rank of $\theta_V$ drops to one on the singular locus of 
$C_4$, 
which must have codimension $6$, hence be a finite set, over which the fibers are projective lines. The fact that $\theta_V$ has a two dimensional kernel means that we can find a basis $v_1,v_2,v_3$ of $V$ such that $\theta_i(v_1,v_2)=\theta_i(v_1,v_3)=0$ for all $i$. 
Completing with two vectors $v_4,v_5$ and taking the dual basis, we conclude that every $\theta_i$ belongs to the space of forms generated 
by  $v_1^\vee \wedge v_4^\vee$, $v_1^\vee\wedge v_5^\vee$ and $\wedge^2(v_1^\perp)$. In particular $\langle \theta_1, \theta_2, \theta_3, \theta_4\rangle$ has to meet   $\wedge^2(v_1^\perp)$ in dimension at least two, which means that $V$ defines a pair of planes $\pi_p, \pi_q$ in $P_0$, whose intersection point is a line in $V$. Finally, $V$ defines a hyperplane $H_{pq}$ of $V_4$, 
and the corresponding fiber is the set $G(2,H_{pq})\simeq \PP^2$ of planes in $H_{pq}$.

Conversely, such a pair of planes being given, we can decompose $\theta$ is an adapted basis as
$$\theta=e_1^\vee\otimes f_1^\vee\wedge f_2^\vee+e_2^\vee\otimes f_3^\vee\wedge f_4^\vee+e_3^\vee\otimes \theta_3 +e_4^\vee\otimes \theta_4,$$
and then the conditions $\theta_3(f_5,\bullet)=\theta_4(f_5,\bullet)=0$ define a $3$-plane $V$ containing $f_5$. This exactly means that 
the singular locus of $C_4$ consists in ten points, in natural bijection with $L_0$. 
\qed

\begin{prop}
Each singular point of $C_4$ defines a plane in the Segre cubic primal $C_3$. The five remaining 
planes are the projectivized kernels of the five singular form $\omega_1, \ldots, \omega_5$. 
\end{prop}

\proof By definition, a point $[v]\in \PP(V_5)$ belongs to $C_3$ when the four linear forms $\theta_i(v,\bullet)$
on $V_5$ are linearly dependent. In the proof above, we have seen that a singular point in $C_4$ corresponds to a 
three-plane $V=\langle v_1,v_2,v_3\rangle$ in $V_5$ with $\theta(v_1,v_2)=\theta(v_1,v_3)=0$. So for any $v\in V$, 
the linear forms $\theta_i(v,\bullet)$ vanish on $v_1$, and also on $v$ by skew-symmetry. When $v$ and $v_1$ are independent, the four linear forms $\theta_i(v,\bullet)$ therefore belong to the three-dimensional space $\langle v,v_1\rangle^\perp\subset V_5^\vee$, so they must be linearly dependent. Hence $\PP(V)\subset C_3$.

That the projectivized kernels $\PP(K_j)$ of the five singular skew forms $\theta_j$ are contained 
in $C_3$ is obvious, since $\theta_j(v,\bullet)=0$ for $v\in K_j$ is a linear dependence relation 
between the $\theta_i(v,\bullet)$.\qed 

\medskip 
Note that we also have a special point $[v_1]$ in each of the ten planes $\PP(V)$. Moreover the five planes 
$\PP(K_1),\ldots, \PP(K_5)$ meet pairwise at a single point. In particular, they provide one of the special 
subcollections of the Cremona-Richmond configuration described in Proposition \ref{subconfigurations}. 




\medskip
Also observe that a form $\omega$ which is as above in the span of $v_1^\vee\wedge v_4^\vee$,
$v_1^\vee\wedge v_5^\vee$ and $\wedge^2(v_1^\perp)$, but does not belong to $\wedge^2(v_1^\perp)$, can be written as $v_1^\vee\wedge w^\vee+\gamma$ with $\gamma\in \wedge^2(v_1^\perp)$ and $w^\vee$ a combination of $v_4^\vee$ and $v_5^\vee$. It has rank two when $\gamma$ has rank (at most) two and $w^\vee\wedge\gamma=0$, which means if $w^\vee\ne 0$ that $\gamma$ is divisible by  $w^\vee$. 
But then $\omega$ itself is divisible by $w^\vee$, and since $w^\vee$ is a combination of $v_4^\vee$ and $v_5^\vee$
this implies that $\omega(v_2,v_3)=0$. In other words, the linear form that defines $H_{pq}\subset V_4$ vanishes at the point that 
corresponds to $\omega$. This exactly means that 
$$p_i\in H_{jk} \qquad i\ne j,k.$$
We thus get in $G(2,V_4)$ a collection of $5+10$ planes, such that each plane of the second type meets exactly three planes of the first type. Hence a configuration $(10_3,5_6)$. The condition that $(jk)$ be disjoint from $(lm)$, so that the two hyperplanes meet in $p_n$, defines a copy of the Petersen graph. 

\medskip\noindent  
Being a degeneracy locus of a morphim between vector bundles, $C_4$ admits two natural resolutions of 
singularities; $X_4$ is one of them. For the other one, we need to impose a rank one kernel in the
source of the morphism $\wedge^2V\ra V_2^\vee$; note that a rank one subspace of $\wedge^2V$ is always  
of the form $\wedge^2W$ for $W\subset V$ a rank two subspace. But then the composition $\wedge^2W\ra V_2^\vee$
vanishes exactly when $W$ defines a point in the del Pezzo surface $X_2\subset G(2,V_5)$. Our second 
resolution of singularities is thus simply $\PP_{X_2}(Q)$, the projectivisation of the quotient bundle 
of $G(2,V_5)$, restricted to $X_2$. The two resolutions are dominated by $\tilde{X}_4$, the set of triples 
$(U_2,V_3\supset W_2)$ such that $(U_2,V_3)$ belongs to $X_4$ and $W_2$ belongs to $X_2$. We get a diagram:

$$\xymatrix{& & \tilde{X}_4 \ar[dl]_\alpha\ar[dd]^\pi\ar[dr]^\beta & & \\
& X_4\ar[dl]_{p_1}\ar[dr]_{p_2} & & \PP_{X_2}(Q)\ar[dl]^{q_2}\ar[dr]^{q_1} &  \\
G(2,V_4) & & C_4 & & X_2. 
}$$

\smallskip

\begin{prop}
The morphism $q_2: \PP_{X_2}(Q)\lra C_4$ is a small resolution of singularities, contracting ten lines to the ten 
singular points of $C_4$. These ten lines are mapped by $q_1$ to the ten lines in the del Pezzo surface $X_2$. 

The morphism $\beta$ is the blow-up of the ten exceptional lines of $q_2$, as well as $\alpha$ is the blow-up 
of the ten exceptional planes of $p_2$. 

Finally, $\pi$ is the blow-up of the ten singular points of $C_4$, its 
exceptional divisor being the disjoint union of ten copies of $\PP^2\times \PP^1$.
\end{prop}

\medskip\noindent {\it Remark.} Contrary to $X_4$, the fourfold $X_4'=\PP_{X_2}(Q)$ is not Fano but only weak Fano.
Indeed, the canonical bundle of $X_2$ is $\det(Q^\vee)$, so the canonical bundle of $X_4'$ is $\cO_{X'_4}(-3)$. 
The quotient bundle $Q$ is obviously not ample on $G(2,V_5)$, and neither is it when restricted to $X_2$ since 
the morphism defined by $\cO_{X'_4}(1)$ is precisely $q_2$ and has non trivial fibers. But of course $Q$ is obviously 
nef, and it is also big since 
$$\int_{X'_4}\cO_{X'_4}(1)^4=\int_{X_2}s_2(Q)=\int_{G(2,V_5)}s_2(Q)\sigma_1^4=2>0.$$

\medskip
Note also the striking similarity with the two main projective resolutions of the Segre cubic, which can be 
encapsulated in a similar diagram
$$
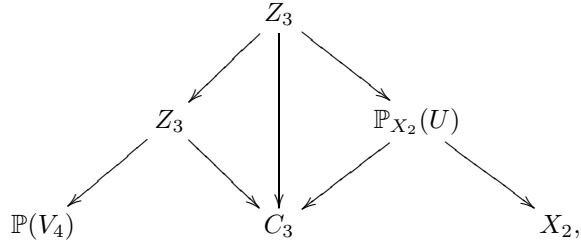
\xymatrix{& & \tilde{Z}_3 \ar[dl]\ar[dd]\ar[dr] & & \\
& Z_3\ar[dl]\ar[dr] & & \PP_{X_2}(U)\ar[dl]\ar[dr] &  \\
\PP(V_4) & & C_3 & & X_2,
}$$

\noindent 
where $Z_3$ is the blowup of $\PP(V_4)=\PP^3$ at five points. Two important differences: $Z_3$, contrary to $X_4$, is 
only weak Fano; $Z_3$ and $ Z'_3=\PP_{X_2}(U)$, contrary to $X_4$ and $X'_4$, are related by flops and therefore
derived-equivalent. Instead of that, we have:

\begin{prop}
The birational map $q_2\circ p_2^{-1}: X_4\dashrightarrow X'_4$ is a flip. 
\end{prop}

\proof Since $X_4$ is Fano, we need only to check that the canonical bundle of $X'_4$ is nef on the non trivial 
fibers of the projection to $C_4$. But we have seen that  $K_{X'_4}=\cO_{X'_4}(-3)$, the fiber of $\cO_{X'_4}(-1)$
at a point defined by a flag $U_2\subset U_3$ being $U_3/U_2$. On a fiber $F$ of the projection to $C_4$, by 
definition $U_3$ is fixed, so $\cO_{X'_4}(-1)_{|F}$ is base point free, hence also $K_{X'_4|F}$. \qed 

\medskip According to the Bondal-Orlov conjecture, there should therefore exist a fully faithful functor 
$D^b(X'_4)\lra D^b(X_4)$ that would be interesting to describe explicitely. 



\subsection{Pencils of skew-forms and the first projection} 

In order to describe the projection to $G(2,V_4)$, we first note that a plane in $V_4$ defines through $\theta$ 
a pencil of skew-symmetric forms in five variables, and that such pencils have been classified. 
In fact, for a two dimensional 
vector space $V_2$, the action of $GL(V_2)\times GL(V_5)$ on $V_2^\vee\otimes \wedge^2V_5^\vee$ has finitely many orbits,
which are described in \cite{kw-E6}. Let us only mention that there are exactly eight orbits:  
the open orbit $\of_7$, 
an orbit $\of_6$ of codimension two and another $\of_5$ of codimension four, and then all the other orbits have bigger codimension. 

The orbit $\of_5$ (or rather its closure) is characterized as consisting of tensors of rank at most four, in the sense that they belong to  $V_2^\vee\otimes \wedge^2V_4$ for some hyperplane $V_4\subset V_5^\vee$. 
The orbit $\of_6$ (or rather its closure) is characterized as consisting of those 
pencils in $\wedge^2V_5^\vee$ admiting a rank two element. 
So the open orbit $\of_7$ parametrizes pencils of forms of constant rank four. By \cite[Proposition 2]{manmez}, given such a pencil one can find a basis of $V_5$ for which the two skew-forms 
$$\omega_1=f_1^\vee\wedge  f_3^\vee+f_2^\vee\wedge  f_4^\vee, \qquad \omega_2=f_1^\vee\wedge  f_4^\vee+f_2^\vee\wedge  f_5^\vee$$
are generators. The projective line $\langle f_1^\vee, f_2^\vee\rangle $ is the {\it pivot} of the pencil. Now, observe that if a three-plane $V\subset V_5$ is isotropic with respect to any skew-form $s\omega_1+t\omega_2$ of the pencil, it has 
to contain its kernel, which is generated by $s^2f_3-stf_4+t^2f_5$. So necessarily $V=\langle f_3,f_4,f_5\rangle$, 
the orthogonal to the pivot. 

\begin{prop}\label{p1}
The projection of $X_4$ to $G(2,V_4)$ is birational. The exceptional 
locus in $G(2,V_4)$ is the union of five planes, intersecting in the ten points of 
 $L_0$, whose fibers are quadratic surfaces.
\end{prop}

\proof 
The fiber of the projection $p_1: X_4\lra G(3,V_5)$
over the point $[U]\in G(2,V_4)$ is defined by the morphism $\theta_U: U\ra \wedge^2V_5^\vee$. 
This morphism is injective and  we thus get a pencil of skew-symmetric forms. If this pencil is generic, 
which means that it has constant rank, then we have just seen that there is a unique three-plane in $V_5$
which is isotropic with respect to any skew-form in the pencil. This three-plane is the image of the induced map
$\theta_U^{(2)}: S^2U\ra \wedge^4V_5^\vee\simeq V_5$. In particular, $p_1$ is birational. 


Special fibers will occur when the pencil $Im(\theta_U)$ becomes special in some way.  By the usual arguments for 
orbital degeneracy locus \cite{ODL1}, we need to take into account, in the space of pencils, only those orbits of codimension smaller than five, which apart from the open orbit are the orbits $\of_5$ and $\of_6$ we have described above. 

Pencils in $\of_5$ contain two skew-forms of rank two. In our case, they must be two of the skew-forms 
$\omega_1,\ldots ,\omega_5$, say $\theta_1$ and $\theta_2$. 
Choose an adapted basis such that $\theta_1=f_1^\vee\wedge f_2^\vee$
and $\theta_2=f_3^\vee\wedge f_4^\vee$, so that
$$\theta_U=e_1^\vee\otimes f_1^\vee\wedge f_2^\vee+e_2^\vee\otimes f_3^\vee\wedge f_4^\vee.$$ 
It is straightforward to check that the three-planes that are isotropic 
with respect to any skew-form in the pencil are those generated by $f_5$, a vector in $\langle f_1,f_2\rangle$, 
and a  vector in $\langle f_3,f_4\rangle$. We thus get for fiber a copy of $\PP^1\times \PP^1$. 

Finally, pencils in  $\of_6$ contain exactly one skew-form of rank two, say $\theta_1$. 
To describe the corresponding fiber we must understand the $3$-planes isotropic with respect to both the generic form 
$\theta_2$ and the degenerate form $\theta_1=f_1^\vee\wedge f_2^\vee$. Such a $3$-plane
must contain the kernel of $\theta_2$; let us choose a generator $f_5$
and a hyperplane $H_4$ in $V_5$ not containing $f_5$. We may suppose that $f_2^\vee$
vanishes on $f_5$. The $3$-planes we are looking for are in correspondence with 
the $2$-planes $H=\langle h, h'\rangle$ in $H_4$ such that $\omega_2(h,h')=0$ and $f_1^\vee(h)=f_1^\vee(h')=0$. 
Such a $2$-plane must be contained in 
the kernel $K_3$ of  $f_1^\vee$, and it has to contain the kernel $K_1$
of the restriction of $\omega_2$ to $K_3$. We finally get for fiber a pencil of planes. 
\qed

\medskip
To summarize, the exceptional locus is the union of five planes $\pi_1,\ldots ,\pi_5$ in $G(2,V_4)$, 
where $\pi_i$ parametrizes the planes in $V_4$ containing $\omega_i$. Any two of these five planes meet
at a single point, over which the fiber of $p_1$ is a quadratic surface. 

If $U_2$ does not belong to any of the five exceptional planes, we have seen that $U_3$ is the span of the 
kernels of the two-forms $\theta(v)$, for $v\in U_2$. Since this kernel can be computed 
as $\theta(v)\wedge\theta(v)$,  there is a natural 
associated conic bundle over $G(2,V_4)$ minus the five exceptional 
planes. This also stresses the analogy with the construction of the Segre primal $C_3$ as the image of a rational map 
$\PP(V_4)\dasharrow \PP(V_5)$ defined by $\theta$. Here we get $C_4$ as the image of  a rational map 
$G(2,V_4)\dasharrow G(2,V_5)$ also defined by $\theta$. We will put its equations in simple form in the next section. 




\subsection{Blow-up and contract}
Proposition \ref{p1} suggests to construct $X_4$
 by first blowing-up $G(2,V_4)$
along the $10$ points of $L_0$, then the strict transforms of the $5$ planes, which are Del Pezzo surfaces of degree five. 
The first blow-up $Bl_0: G_0\lra G(2,V_4)$
gives $10$  exceptional divisors $E_{ij}\simeq \PP^3$ for $1\le i<j\le 5$, each with a pair 
of skew lines $\ell_i,\ell_j$ coming from the two planes $\pi_i$ and $\pi_j$ intersecting at $p_{ij}$. 
The second blow-up $Bl_P: G_1\lra G_0$
produces five other exceptional divisors $F_k$ for $1\le k\le 5$, while the strict transform of $\tilde{E}_{ij}$ 
of $E_{ij}$ is the blowup of $E_{ij}$ along $\ell_i\cup\ell_j$. Since the blowup of $\PP^3$ along two skew lines 
is the total space of $\PP(\cO(-1,0)\oplus \cO(0,-1))$ over $\PP^1\times \PP^1$, we deduce that the rational
map to $X_4$ is a morphism. More precisely, it has to coincide with the blowup $Bl_Q: G_1\lra X_4$ 
of the ten quadratic surfaces $S_{ij}=p_1^{-1}(p_{ij})$ in $X_4$. This  explains in particular why the 
Picard number is equal to $6$. 

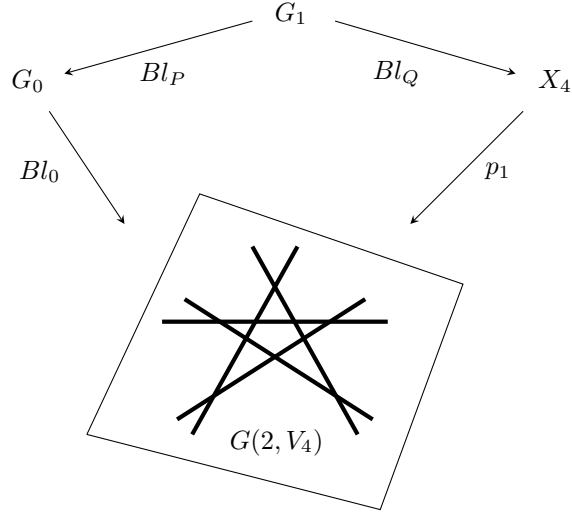
\begin{center}
    \begin{tikzpicture}
      \draw [thin] (.5,-1) -- (2,2.2) -- (5.5,1) -- (4.4,-2) -- (.5,-1);
      \draw [ultra thick] (1.5,.5) -- (4.5,.5); 
         \draw [ultra thick] (3.3,1.5) -- (1.9,-1); 
                \draw [ultra thick] (2.7,1.5) -- (4.1,-1); 
    \draw [ultra thick] (1.7,-.8) -- (4.2,.8); 
     \draw [ultra thick] (4.3,-.8) -- (1.8,.8); 
       \node[text width=4mm] at (2.6,-1.1) {$G(2,V_4)$};
         \draw [stealth-](4.8,1.8) -- (6.3,3.3);
          \node[text width=4mm] at (6.7,3.7) {$X_4$};
          
         \draw [stealth-](1,1.8) -- (0,3.3);
          \node[text width=4mm] at (-.3,3.7) {$G_0$};
          \node[text width=4mm] at (-.2,2.5) {$Bl_0$};
          
          \node[text width=4mm] at (3.2,4.6) {$G_1$};
          
          \node[text width=4mm] at (1.4,3.8) {$Bl_P$};
          \node[text width=4mm] at (4.5,3.8) {$Bl_Q$};
          
          \node[text width=4mm] at (6,2.5) {$p_1$};
         \draw [-stealth](2.7,4.5) -- (.2,3.8);
         \draw [-stealth](3.8,4.5) -- (6.2,3.8);
    \end{tikzpicture}
\end{center}


\medskip 
Let $F=F_1+\cdots +F_5$, and let $E$ be the sum of the ten divisors $\tilde{E}_{ij}$ in $G_1$. 
from the identity 
$$K_{G_1}=-4H_1+3E+F=K_{X_4}+E=-H_1-H_2+E$$
we deduce the relation $3H_1=H_2+2E+F$. 

\medskip The exceptional locus of $p_2$ defines  a collection of $10$ planes in $X_4$, 
contracted to the ten singular points of $C_4$, and that we can identify with their isomorphic images in $G(2,V_4)$.
Recall that in this Grassmannian we have the five planes $\pi_1,\ldots ,\pi_5.$

\begin{prop} The resulting collection of $10+5$ planes in $G(2,V_4)$ is in natural correspondance with 
the Cremona-Richmond configuration.
\end{prop} 

\subsection{Incidences with the Segre cubic}
Now we relate the two varieties $X_3$ and $X_4$ by considering the 
incidence correspondence 
$$I=\{(A_1,B_1),(U_2,U_3)\in X_3\times X_4, \; A_1\subset U_2, \; 
B_1\subset U_3\}.$$
Recall that by definition, $B_1$ is (contained in) the kernel of 
$\theta(v)$ for $v\in A_1$, while $U_3$ is the linear span of the kernels
of the two forms $\theta(u)$ for $u\in U_2$; this kernel depends quadratically on $u$ since it is given by $\theta(u)\wedge\theta(u)$.  
This implies that $I$ is (generically) a $\PP^1$-bundle 
over $X_4$, 
and (generically) a $\PP^2$-bundle over $X_3$. We have a commutative diagram 
$$\xymatrix{& I\ar[dl]\ar[d]\ar[dr] & \\
X_3\ar[d] & Fl(1,2,V_4)\ar[dl]\ar[dr] & X_4\ar[d] \\
\PP(V_4) & & G(2,V_4). 
}$$
\medskip 

One easily checks that:

\begin{prop} 
The base locus of the birational projection $I\lra Fl(1,2,V_4)$ is the union of five disjoint planes.
\end{prop}

\proof Denote by $\psi_i$ the plane $G(2,V_4)$ parametrizing the lines
in $\PP(V_4)$ that pass through $p_i$, and by $\Psi_i$ its lift in 
$ Fl(1,2,V_4)$. The preimage in $I$ of a point in $\psi_i$ is given by a flag 
$B_1\subset U_3\subset V_5$ such that $U_1$ is contained in the 
kernel $K_i$ of $\omega_i$ and $U_3$ contains




\subsection{Projective duality} 

We have seen that $X_4$ is birationally equivalent to the projective bundle $\PP_{X_2}(Q)$ over the del 
Pezzo surface $X_2$. Since $Q^\vee$ has no section, we would rather write it as $\PP=\PP_{X_2}(\wedge^2Q^\vee)$, 
in which case the relative tautological bundle $\cO_\PP(-1)$ sends $\PP$ to $\PP(\wedge^2V_5^\vee)\simeq 
\PP(\wedge^3V_5)$, the image being $C_4\subset G(3,V_5)$. We are then in the context of Homological 
Projective Duality for projective bundles, according to which $\PP\ra\PP(\wedge^2V_5^\vee)$ is dual to 
 $\PP^*\ra\PP(\wedge^2V_5)$, with $\PP^*$ the projective bundle $\PP_{X_2}(W\wedge V_5)$, where $W$ denotes
 the rank to tautological bundle. 
 
 \begin{prop} 
 The image of $\PP^*\ra\PP(\wedge^2V_5)$ is an octic hypersurface in $ \PP(\wedge^2V_5)$, 
 containing the Grassmannian $G(2,V_5)$ in its singular locus. 
 \end{prop}

\proof First consider the full projective bundle $\PP_{G(2,V_5)}(W\wedge V_5)$ and its projection to 
$\PP(\wedge^2V_5)$. The generic fiber is a copy of $\QQ^3$ (while the special fibers, that occur over $G(2,V_5)$,
are codimension two Schubert cycles). When we restrict to $X_2$, we cut the fibers by linear spaces of codimension four. Generically, they meet the span of the fiber at one point; in codimension one, this point will be on the fiber itself.
This implies that $\PP^*\ra\PP(\wedge^2V_5)$ is birational onto its image, which must be a hypersurface. As usual,
we compute the degree of this hypersurface as 
$$\int_{\PP^*}\cO_{\PP^*}(1)^8=\int_{X_2}s_2(W\wedge V_5)=\int_{G(2,V_5)}(2\sigma_2+\sigma_{11})\sigma_1^4=8.$$
Here we used exact sequences to compute the Segre class $s(W\wedge V_5)=c(Q)^5c(S^2U)$, with $c(Q)=1+\sigma_1+\sigma_2+\sigma_3$ and $c(U)=1-\sigma_1+\sigma_{11}$. 

Over a point $W^0$ of the Grassmannian, the fiber of $\PP_{G(2,V_5)}(W\wedge V_5)\ra\PP(\wedge^2V_5)$
is the Schubert cycle of planes $W$ meeting $W^0$ along at least a line. It is desingularized by a $\PP^3$-bundle 
over $\PP(W^0)$. If we fix a line $L\subset W^0$, there exists a plane $W\supset L$ in $X_2$ if and only if 
the four linear forms $\theta(L,\bullet)$ on $V_5/L$ are linearly dependent. This defines a section of 
$\wedge^4(Q^\vee(1))=\cO(3)$ over $\PP^1$, and we conclude that the general fiber of $\PP_{X_2}(W\wedge V_5)\ra\PP(\wedge^2V_5)$ over $G(2,V_5)$ consists in three points. Since this morphism is birational onto
its image, Zariski's main theorem implies that $G(2,V_5)$ is contained in the singular locus. \qed 


\section{Symmetries} 

The symmetries of the Segre cubic primal must be reflected in $X_4$. In this section we describe the
symmetries of $X_4$ in some detail. In particular we will prove:

\begin{prop}
The generic stabilizer of the action of  $PGL(V_4)\times PGL(V_5)$ on $\PP(V_4^\vee\otimes\wedge^2V_5^\vee)$ is the symmetric group $\mathcal{S}_5$.
\end{prop}

What is classically known, as we mentionned in the introduction, is that the action of  $PGL(V_4)\times PGL(V_5)$  on $V_4^\vee\otimes\wedge^2V_5^\vee$ is prehomogeneous. 
The representative of the open orbit given in \cite{ozeki79} is
$$\theta=e_1^\vee\otimes (f_{25}-f_{34})+e_2^\vee\otimes (f_{15}-f_{24})+e_3^\vee\otimes (f_{23}-f_{14})+e_4^\vee\otimes (f_{45}-f_{12}),$$
with the notation $f_{ij}=f_i^\vee\wedge f_j^\vee$. 
The corresponding points in $\PP(V_4)$ and rank two forms are easy to identify; we get
$$\begin{array}{lll}
p_1=e_2+ie_4, & \qquad & \omega_1= (f_1+if_4)\wedge (f_2+if_5), \\
p_2=e_2-ie_4, & \qquad & \omega_2= (f_1-if_4)\wedge (f_2-if_5), \\
p_3=e_1+e_3+e_4, & \qquad & \omega_3= (f_2+f_4)\wedge (f_1+f_3+f_5), \\
p_4=e_1+je_3+j^2e_4, & \qquad & \omega_4= (f_2+j^2f_4)\wedge (f_1+j^2f_3+jf_5), \\
p_5=e_1+j^2e_3+je_4, & \qquad & \omega_5= (f_2+jf_4)\wedge (f_1+jf_3+j^2f_5).
\end{array}$$
Here $j$ and $i$ are primitive fourth and third roots of unity. Each pair $\omega_p,\omega_q$ defines two planes in 
$V_5^\vee$ whose common orthogonal is a line $[e_{pq}].$ Then the planes of the Cremona-Richmond configuration are 
obtained as follows: $P_{pq}$ is generated by the three points $e_{ij}, e_{jk}, e_{ik}$ for $ijk$ distinct from $pq$; 
and $P_p$ is generated by the four points $e_{ip}$ for $i\ne p$. Explicitely, the ten vectors $e_{pq}$ 
can be chosen as follows:
$$\begin{array}{lcl}
e_{12}=(0,0,1,0,0) & \qquad & e_{24}=(i,-j^2,-2ij,1,ij^2) \\

e_{13}=(1,-i,-2,i,1) & \qquad & e_{25}=(i,-j,-2ij^2,1,ij)  \\

e_{14}= (-i,-j^2,2ij,1,-ij^2) & \qquad & e_{34}= (1,0,j^2,0,j)\\

e_{15}=(-i,-j,2ij^2,1,-ij)& \qquad & e_{35}=(1,0,j,0,j^2) \\

e_{23}= (1,i,-2,-i,1)& \qquad & e_{45}=(1,0,1,0,1). 
\end{array}$$

Each $\omega_i$ defines a plane $\pi_i$ in $V_5^\vee$, from which we can deduce a collection of hyperplanes  
$\pi_{ij}=\pi_i+\pi_j$ and points $p_{ijk}=\pi_i\cap (\pi_j+\pi_k).$

\begin{prop} 
For any permutation $i,j,k,l,m$ of $1,\ldots, 5$, $p_{ijk}=p_{ilm}$.
\end{prop}

\proof Explicit check.\qed 

\medskip
We have no convincing explanation of this coincidence, but 
as a consequence, we don't get thirty but only fifteen points in $\PP(V_5^\vee)$. Obviously, $p_{ijk}$ belongs to $\pi_i$, 
hence to any of the four hyperplanes $\pi_{il}$, $l\ne i$. Conversely, $\pi_{ij}$ contains the three points $p_{iab}$ plus the three 
points $p_{jcd}$.  

\begin{prop} 
The fifteen points $\pi_{ijk}$ and the ten hyperplanes $\pi_{ij}$ in $\PP(V_4)$ 
form a configuration $(15_4,10_6)$. 
\end{prop} 

We thus recover the abstract configuration classically defined by the Segre primal.
In particular the  fifteen points $\pi_{ijk}$ should be in natural correspondence with planes in the Segre primal. 

\medskip
Automorphisms in $PGL(V_4)\times PGL(V_5)$ that fix $\langle \theta\rangle$ are in bijective correspondence with 
elements of $PGL(V_5)$ fixing the four-plane generated by the $\omega_i$'s. Automatically such an automorphism
will preserve  the set of five planes $\pi_1,\ldots ,\pi_5$, hence the collection of the thirty points $p_{ijk}$.

In order to show that any permutation of the five planes can be lifted to $PGL(V_5)$, it is enough to lift two
generators of $\mathcal{S}_5$, say a transposition and a complete cycle. 
By sending $f_i$ to $\epsilon_if_i$ with $\epsilon_i=1$ for $i$ odd and  $\epsilon_i=-1$ for $i$ even, we exchange $\pi_1$ and $\pi_2$ and let the 
three other planes fixed. So let us turn to a maximal cycle. We claim that the cycle $(12345)\in\mathcal{S}_5$ can be 
lifted to the tranformation of $GL(V_5)$ given by 
$$\begin{array}{rcl} 
f_1 & \mapsto &   \frac{j}{3}f_1-2ijf_2+ \frac{j}{3}f_3-ijf_4+ \frac{4j}{3}f_5,\\
f_2 & \mapsto &  -\frac{2i}{3}f_1-f_2+\frac{i}{3}f_3+\frac{i}{3}f_5, \\
f_3 & \mapsto & \frac{4j^2}{3}f_1+4ij^2f_2-\frac{2j^2}{3}f_3-4ij^2f_4+\frac{4j^2}{3}f_5, \\
f_4 & \mapsto & -\frac{ij}{3}f_1-\frac{ij}{3}f_3-jf_4+\frac{2ij}{3}f_5,\\
f_5 & \mapsto & \frac{4}{3}f_1+if_2+\frac{1}{3}f_3+2if_4+\frac{1}{3}f_5.
\end{array}$$

\smallskip

\begin{coro}
The automorphism group of the Fano fourfold $X_4$ is $Aut(X_4)=\mathcal{S}_5$.
\end{coro}

\proof An automorphism of $X_4$ is induced by a linear transformation in $PGL(V_4)\times PGL(V_5)$ preserving $\theta$. 
Considered as a homomorphism from $V_4$ to $\wedge^2V_5^\vee$, $\theta$ defines a codimension four linear section of 
$G(2,V_5)$, that is a degree five del Pezzo surface $S_5$. This implies that $Stab(\theta)$ embeds into $Aut(S_5)$,
which is well-known to be $\mathcal{S}_5$. Since we know by the previous computations that $Stab(\theta)$ contains 
$\mathcal{S}_5$, we are done.
\qed 

\medskip Once we identify $\mathcal{S}_5$ with the stabilizer of $\theta$ in $SL(V_4)\times SL(V_5)$, we get actions of 
$\mathcal{S}_5$ on $V_4$ and $V_5$, clearly irreducible. Up to the sign representation there is a unique irreducible
representation of  $\mathcal{S}_5$  of dimension $4$, and a unique one of dimension  $5$. 
The complex (\ref{lascoux}) shows that $C_4$ is cut out by 
a family of quadrics on $G(2,V_5)$ parametrized by $V_4$, hence a $\mathcal{S}_5$-invariant copy of $V_4$ inside 
$S_{22}V_5$. We will show later on that this copy  is unique.
(This point of view from finite group representation theory is typically used in \cite{dolgachev18}.
Something with the same flavour has been done 
in \cite{bc} for the quintic del Pezzo surface.)

\medskip
We use the character table of  $\mathcal{S}_5$ (see for example \cite{dfl}) to compute some plethysm and tensor product representations. Recall that  $\mathcal{S}_5$ has irreducible representations of dimension $1,1,4,4,5,5,6$ that we 
denote by $U_1, U_1^-, U_4, U_4^-, U_5, U_5^-, U_6$. All these representations are self-dual, being defined over the real numbers.
Concretely, $U_1$ is the trivial representation,  $U_1^-$ is 
the sign representation, $U_4$ is the natural representation, $U_4^-=U_4\otimes U_1^-$ and $U_6=\wedge^2U_4$. 
One computes that 
$$S^2U_4=U_5\oplus U_4\oplus U_1, \qquad
\wedge^2U_5=U_4^-\oplus U_6.$$
The last decomposition implies in 
particular that $(U_4^-)^\vee\otimes\wedge^2U_5^\vee$ contains a unique $\mathcal{S}_5$-invariant tensor $\theta_{\mathcal{S}_5}$, up to scalars. 

At this point it could therefore be reasonable to reverse the whole process and start from the representation
theory of $\mathcal{S}_5$. One should be able to check directly that  $\theta_{\mathcal{S}_5}$ is generic, and then we should get  $\theta_{\mathcal{S}_5}$-invariant descriptions of all the objects we have been studying. 

Note that $S^2U_4=U_5\oplus U_4\oplus U_1$ allows to construct $U_5$ from $U_4$, as the space of quadrics 
which are apolar to the obvious invariant cubic. In coordinates $x_1,\ldots, x_5$ permuted by 
 $\mathcal{S}_5$, the representation $U_4$ is the hyperplane 
 $x_1+\cdots + x_5=0$, the invariant cubic is  $x_1^3+\cdots + x_5^3$ 
 and the apolar quadrics are of the  form $\sum_{i\ne j}a_{ij}x_ix_j$ with
 $$a_{ij}=a_{ji} \; \forall i\ne j, \qquad \qquad  \sum_{i\ne k}a_{ik}=0 \;\forall k.$$
 We get ten indeterminates and five independent relations, consistently with the fact 
 that these quadrics should span a copy of $V_5$. 

Inside the space $V_5$ of apolar quadrics to the invariant cubic, note that
we have  $q_{ij,kl}=(x_i-x_j)(x_k-x_l)$ for $i,j,k,l$ distinct integers. 
These quadrics  are subject to the Plücker type relations $q_{ij,kl}-q_{ik,jl}+q_{il,jk}=0$. 
This suggests to define the following elements of $\wedge^2V_5$:
$$\begin{array}{rcl} 
Q_1 & = & q_{23,45}\wedge q_{24,35}, \\
Q_2 & = & q_{13,45}\wedge q_{14,53}, \\
Q_3 & = & q_{12,45}\wedge q_{14,25}, \\
Q_4 & = & q_{12,35}\wedge q_{13,52}, \\
Q_5 & = & q_{12,34}\wedge q_{13,24}.
\end{array}$$
Obviously, for any permutation $\sigma \in \mathcal{S}_5$ one must have $\sigma(Q_i)=\pm Q_{\sigma(i)}$. 
We also let, for a pair $i\ne j$ with complement $p,q,r$ in $1\ldots 5$, 
$$Q_{i,j}=q_{ip,qr}\wedge q_{jp,qr}+q_{iq,rp}\wedge q_{jq,rp}+q_{ir,pq}\wedge q_{jr,pq}. $$

\begin{prop}
The action of $\cS_5$ on $\langle Q_1,\ldots , Q_5\rangle$ gives a copy of the representation $U_4^-$ in $\wedge^2U_5$. 
Similarly, the action of $\cS_5$ on $\langle Q_{i,j}, 1\le i< j\le 5\rangle$ gives a copy of the representation $U_6$.
\end{prop}

What have we gained in doing all that? 
First, we get a better, more symmetric normal form for the generic $\theta$ than that of Ozeki, as
$$\theta_{\mathcal{S}_5}=e_1\otimes Q_1+e_2\otimes Q_2+e_3\otimes Q_3+e_4\otimes Q_4+e_5\otimes Q_1,$$
with $e_1+\cdots +e_5=0$. 

\smallskip
Also, we can make explicit the quadratic equations of $\cC_4$. A character computation yields:

\begin{lemma}
The multiplicity of $U_4^-$ inside $S^2(\wedge^2U_5)$ is equal to one. 
\end{lemma}

So the space of quadratic equations we are looking for is uniquely defined in terms of the $\cS_5$-action. 
Moreover, recall that $\wedge^2U_5=U_4^-\oplus U_6$. Another character computation shows that the copy 
of $U_4^-$ that we are looking for inside $S^2(\wedge^2U_5)$ is in fact contained inside 
$U_4^-\otimes U_6=U_4^-\otimes \wedge^2(U_4^-)\subset U_4^-\otimes End(U_4^-)$ (recall that $U_4^-$ is self-dual).
So there is an obvious map to  $U_4^-$, and dually, this says that the space of quadrics we are looking for 
is generated by the five quadrics 
$$CQ_i=\sum_{j\ne i} Q_{i,j}Q_j, \qquad 1\le i\le 5.$$

\medskip\noindent {\it Remark}.
Since $Aut(C_3)= \mathcal{S}_6$,  certain automorphisms of the Segre primal do not lift to $X_4$.
would it be possible that $\mathcal{S}_6$ act on $X_4$ by birational transformations?

\section{The Chow ring of $X_4$}

In this section we completely determine the Chow ring of $X_4$, with its structure of $\cS_5$-module.
Let us start with the Picard group. 

From the relation $3H_1=H_2+2E+F$ that we found on $G_1$,  we compute that
$$H_1^4=2, \quad H_1^3H_2=6, \quad H_1^2H_2^2=13, 
\quad H_1H_2^3=14, \quad H_2^4=12.$$
The Picard group is generated by $H_1$, $H_2$ and the five components of $F$, which are permuted by 
$\cS_5$. We deduce:

\begin{prop}
The Chow ring of $X_4$ is generated by $A^*(G)$ and the five divisors $F_1, \ldots , F_5$. 
As a representation of $\cS_5$, the Picard group decomposes as  $$Pic(X_4)\otimes_\ZZ \CC=2U_0\oplus U_4.$$
\end{prop}

We know by Proposition \ref{invariants}
that the middle dimensional Chow group $A^2(X_4)$ has dimension $17$, and we expect that 
the invariant part has dimension four, with two classes coming from $G(2,V_4)$ and two other classes 
from $G(2,V_5)$. We will show that are all come (at least over $\QQ$) from products of divisor classes.



\medskip
We compute the multiplicative structure of the Chow ring by 
embedding it in the Chow ring of $G_1$, that we shall now describe. First, 
the Chow ring of $G_0$ is generated by the Chow ring of $G=G(2,V_4)$ and the ten exceptional divisors $E_{pq}^0$
of the blow-up $b_0=Bl_0$, 
such that 
$$(E_{pq}^0)^4=-1, \quad E_{pq}^0E_{p'q'}^0=0 \; \mathrm{for}\; \{p,q\}\ne \{p',q'\}, \quad E_{pq}^0.b_0^*C=0$$
for any class $C\in A^*(G)$ of positive degree. After this first blow-up, the five planes $\pi_1,\ldots ,\pi_5$ give five disjoint 
surfaces $\Sigma_1, \ldots ,\Sigma_5$, each one being a plane blow-up in five points, that is a del Pezzo surface of degree $5$. 
We denote the four exceptional lines in $\Sigma_p$ by $\ell_p^q$, whose image in $G$ is the point $\pi_{pq}$, for $q\ne p$. 

The second blow-up $b_1=Bl_P$ is the blow-up of these five surfaces. We denote by $F_p^1$ the five exceptional divisors, and by $E_{pq}^1$ the strict transforms of the divisors $E_{pq}^0$. 
Since $F_p^1=\PP(N_p)$, for $N_p$ the normal bundle of $\Sigma_p$ inside $G_0$, we need to describe this normal bundle. Recall that when one blows up one point in a smooth 
variety $X$, creating an exceptional divisor $E$ inside the blow-up $Y\stackrel{\pi}{\ra} X$, the tangent exact sequence is $0\ra TY\ra \pi^*TX\ra i_*TE\ra 0$, where $i:E\ra Y$ denotes the inclusion. 
Since the normal bundle of $\pi_p$ inside the Grassmannian $G$ is the quotient bundle $Q$, 
we get the following diagram:

$$\begin{CD}
 @. 0 @. 0 @. 0 @. \\
 @. @VVV @VVV @VVV \\
0 @>>>  T\Sigma_p @>>> TG_{0|\Sigma_p} @>>> N_p @>>> 0 \\
   @. @VVV        @VVV @VVV \\
0 @>>> b_0^*T\pi_p @>>> b_0^*TG_{|\pi_p} @>>> b_0^*Q @>>> 0 \\
 @.  @VVV        @VVV @VVV \\
 0 @>>> \oplus_{q\ne p}T\ell_p^q  @>>> \oplus_{q\ne p}TE_{pq}^0
 @>>> \oplus_{q\ne p} N_p^q @>>> 0 \\
  @.  @VVV  @VVV        @VVV\\
   @. 0   @. 0     @.   0\\
\end{CD}$$

\smallskip
Here we denoted by $N_p^q$ the normal bundle to $\ell_p^q\simeq \PP^1$ inside $E_{pq}^0\simeq \PP^3$, which is just $\cO_{\ell_p^q}(1)\oplus \cO_{\ell_p^q}(1)$. We deduce the Segre class 
$$s(N_p)=s(b_0^*Q)\prod_{q\ne p}c(\cO_{\ell_p^q}(1))^2\in A^*(\Sigma_p).$$
One the one hand, the Segre class $s(Q)$ equals the Chern class of the tautological bundle on $G$, that is $s(Q)=1-H_1+\sigma_{11}$, and the Schubert class $\sigma_{11}$ restricts to zero on $\pi_p$. On the other hand, on the del Pezzo surface $\Sigma_p$ we have  $\cO_{\ell_p^q}(1)=\cO(-\ell_p^q)_{|\ell_p^q}$, from which we get 
the Segre class $s(\cO_{\ell_p^q}(1))=1+\ell_p^q+2(\ell_p^q)^2$. 
Finally, 
$$s(N_p)=1-H_1+2\sum_{q\ne p}\ell_p^q+2\sum_{q\ne p}(\ell_p^q)^2.$$
We can deduce several intersection numbers, since for any class $C_{3-k}$ of degree $3-k$ on $G_0$, we have the classical formulas 
$$F_p^{k+1}b_1^*C_{3-k}=\int_{F_p}F_p^{k}b_1^*C_{3-k}=
(-1)^k\int_{\Sigma_p}s_{k-1}(N_p)C_{3-k}.$$

\begin{lemma}\label{int1}
$$(F^1_p)^4= 8, \quad (F^1_p)^3H_1=-1, \quad (F^1_p)^2H_1^2=-1, \quad F^1_pH_1^3=0.$$
\end{lemma}

Note also that $F^1_p$ does not meet $E^1_{rq}$ for $r,q\ne p$, but it meets $E^1_{pq}$ transverselly along the surface $S_p^q=b_1^{-1}(\ell_p^q)$. Therefore 
$$\cO_{G_1}(E^1_{pq|F^1_p})=
\cO_{F^1_p}(S_p^q)=b^*_1\cO_{\Sigma_p}(\ell_p^q).$$
Applying the previous formula to $C_{3-k}=(E_{pq}^0)^{3-k}$ we get: 

\begin{lemma}\label{int2} $F^1_pE^1_{rq}=0$ if $r,q\ne p$, but 
$$(F^1_p)^3E_{pq}^1=-2, \quad (F^1_p)^2(E_{pq}^1)^2= 1,\quad F^1_p(E_{pq}^1)^3=0.$$
\end{lemma}

\smallskip
On the other hand, $E_{pq}^0$ gets blown-up along the two-skew
lines $\ell_p^q$ and $\ell_q^p$, and its strict transform 
$E^1_{pq}$ is contracted to the quadratic surface $\ell_p^q\times\ell_q^p$ in $X_4$. This surface is also the 
intersection of $F_p$ and $F_q$ in $X_4$, in particular it 
is contained in $F_p$. We deduce, denoting $Bl_Q$ by $c$, that 
$$c^*F_p=F_p^1+\sum_{q\ne p}E^1_{pq}.$$
Summing up over $p$, we get the relation $c^*F=F^1+2E^1$. 

\begin{coro}
$C_4\subset G(3,V_5)$ is the image of $G=G(2,V_4)$ by the 
linear system $|I_\pi(3H_1)|$ of cubics vanishing along the 
union $\pi$ of the five planes $\pi_1,\ldots ,\pi_5$. 
\end{coro}

We have enough information to describe the full 
intersection product on $X_4$.

\begin{prop}\label{intX4} The nonzero intersection numbers among the divisor 
classes $H_1, F_1,\ldots , F_5$ are the following, for $1\le p\ne q\le 5$:
$$F_p^4=12, \quad F_p^3F_q=-2, \quad F_p^2F_q^2=1,  \quad F_p^3H_1 = -1, \quad F_p^2H_1^2=-1,
\quad H_1^4=2.$$
Moreover we always have $H_1F_pF_q=0$ for $p\ne q$ and $F_pF_qF_r=0$ for $p\ne q\ne r\ne p$. 
\end{prop}

\proof The values of $F_p^3H_1$ and $F_p^2H_1^2$ can be computed directly by restricting to a 
general hyperplane or a general 
codimension two section of $G$; then we avoid the points $\pi_{qr}$ and we are reduced to compute the self-intersection 
of the exceptional divisor for the blow-up of a line in a three-dimensional quadric, or a point in a surface. 
Then we can deduce the value of $F_p^4$ by computing the self-intersection of $H_2=3H_1-F$,  which we know is equal to 
$$12=81H_1^2-108H_1^3F+54H_1^2F^2-12H_1F^3+F^4=162-270+60+F^4.$$
This gives $F^4=F_1^4+\cdots +F_5^4=60$, hence $F_p^4=12$. (But note that this is not equal to $(c^*F_p)^4=-4$, as a consequence of the fact that $F_p$ contains four of the quadratic surfaces blown-up by $c$.)

The other intersection numbers can be computed by pulling-back by $c$ and using Lemma \ref{int1}.\qed

\begin{prop}\label{square} 
The square map $S^2A_1(X_4)\lra A^2(X_4)$ is surjective.
As a consequence, the $\cS_5$-module structure  of $A^2(X_4)$ is
$$A^2(X_4)=4U_0\oplus 2U_4\oplus U_5.$$
\end{prop} 

\proof The decomposition of the $\cS_5$-module $S^2A^1(X_4)$ is
$4U_0\oplus 3U_4\oplus U_5$, the sum of three isotypic components, and the kernel of the square map must decompose accordingly. 

First consider the four invariant classes 
$H_1^2, H_1F, F^{(2)},F^{(11)}$, where 
$$F^{(2)}=\sum_p F_p^2, \qquad F^{(11)}=\sum_{p<q} F_pF_q.$$
Suppose that there is a relation $aH_1^2+bH_1F+cF^{(2)}+dF^{(11)}=0$. Multiplying successively
by $H_1^2, H_1F_p, F_p^2, F_pF_q$ and using the results of Proposition \ref{intX4}, we  deduce that $2a-5c=0$, $b+c=0$, $a+b+16c-8d=0$,  
$4c-d=0$, hence $a=b=c=d=0$. 

Now consider the possibility that $U_5$ be contained in the kernel of the square map. We claim that $U_5$ is embedded inside $S^2A^1(X_4)$ as the space of linear combinations 
$\sum_{p\ne q}a_{pq}F_pF_q$ with $a_{pq}=a_{qp}$ and $\sum_r a_{pr}=0$ for all $p,q$. Indeed, this defines an invariant 
five-dimensional subspace of $S^2A^1(X_4)$, not containing any 
invariant class, so it must be $U_5$. A typical element is 
$$3F_pF_q-(F_p+F_q)\sum_{r\ne p,q}F_r+\sum_{s,t\ne p,q}F_sF_t.$$
If this was zero in $A^2(X_4)$, multiplying by $F_pF_q$ would 
imply that the intersection number $F_p^2F_q^2=0$, which is not the case. 

We can conclude that the kernel of the square map must be contained in the isotypical component $3U_4$ of $S^2A^1(X_4)$, 
which is generated by the three copies of $U_4$ respectively obtained as the linear combinations $\sum_p a_pH_1F_p$,  $\sum_p a_pFF_p$ and $\sum_p a_pF^2_p$ for $\sum_p a_p=0$. A copy of $U_4$ in the kernel corresponds to a relation of the form
$$uH_1F_p+vFF_p+wF_p^2=I\qquad \forall p,$$
for $I$ an invariant class. Since $I$ is invariant, 
multiplying by $H_1F_p$ and $H_1F_q$ must then give the same intersection number, which gives the relation $-u-v-w=0$. Similarly, multiplying by $F_p^2$ or $F_q^2$
must give the same result, that is $-u+4v+12w=-v+w$. Finally, multiplying by $F_pF_q$ or $F_qF_r$ with $q,r$ distinct from $p$
must also give the same result, that is $-v-2w=0$. These three equations are linearly dependent and reduce to $u=w$ and $v=-2w$,
which proves that there is a unique copy of $U_4$ in the kernel of the square map. This concludes the proof. \qed 



\subsection*{Threefolds} Consider the two families of divisors in $X_4$ given by sections of $H_1$ and $H_2$,
respectively. Since a general  hyperplane section in $G(2,V_4)\simeq\QQ^4$ will avoid the ten points $\pi_{pq}$, the first 
ones are just blowups of five disjoint lines in $\QQ^3$. For the same reason, the second ones, say $Z_3$, are isomorphic 
with their images in $G(3,V_5)\cap H_2$, which are codimension two degeneraci loci defined by the condition that 
the morphism $\wedge^2V\ra V_4^\vee$ has rank exactly two. Its image is then the pullback from $G(2,V_4)$ of the dual
quotient bundle $Q^\vee$. In particular we get an exact sequence 
$$0\ra \cO(1,-2)\ra p_2^*(\wedge^2V)\ra p_1^*Q^\vee\ra 0$$
\noindent
on $Z_3$. This shows in particular that $\cO(-1,2)$, 
the restriction of $2H_2-H_1=5H_1-2F$, is generated by sections on $Z_3$.
The image in $Z_3$ is the closure of the planes $U\subset V_4$ such that the image of $S^2U\ra\wedge^2V_5\simeq V_5$
is isotropic with respect to some three-form  on $V_5$. This defines a section of $\wedge^3(S^2U)^\vee=
\det(U^\vee)^3$, so that the image of $Z_3$ in $G(2,V_4)$ is a singular cubic hypersurface. 

\subsection*{K3 surfaces} By taking sections of $H_1\oplus H_2$ in $X_4$, we get a family of smooth K3 surfaces $S$ in $X_4$. 
We denote by $h_1$, $h_2$, $f_1,\ldots ,f_5$ the restriction to $S$ of the divisors $H_1, H_2, F_1,\ldots  F_5$. 

\begin{prop} The intersection numbers of these divisors in  $S$ are 
$$h_1^2=6, \quad  h_2^2=14, \quad h_1h_2=13, \quad h_1f_i=1,\quad h_2f_i=5, \quad f_if_j=-2\delta_{ij}.$$
\end{prop}

\proof This is an immediate consequence of the computations above, since for two divisors $A,B$ on $X$ restricting to $a,b$ on $S$, we have  $ab=ABH_1H_2$. \qed 

\medskip An obvious consequence is that $h_1,f_1,\ldots ,f_5$ are linearly independent. Moreover, the curves 
$C_i=F_i\cap S$ are $(-2)$-curves on $S$, mapping to lines on $G(2,V_4)$ and to rational quintics in $G(2,V_5)$. 
The divisor $5h_1-h_2=2h_1+f$ should contract these five $(-2)$-curves to the five singular points 
of a surface $\bar{S}$. Note that this is a divisor of degree $34$, so $\bar{S}$ could be a degeneration of 
a smooth K3 surface of genus $18$. Mukai described the generic such K3 surface as the zero locus in $OG(3,9)$ 
of five sections of the rank two spinor bundle. What is the connection? Note that we have a family of surfaces of 
dimension $5+9=14=19-5$, which is coherent with the expectation that imposing $5$ nodes on a K3 surface of genus $18$
should give five independent conditions. 

\section{The Igusa quartic and the Coble fourfold} 

Given a linear form $h$ on $V_5$, there is an associated quadratic form $Q_h$ on $V_4$:
$$Q_h(v)=h\wedge \theta(v)\wedge \theta(v)\in\wedge^5V_5^\vee\simeq\CC.$$

\begin{prop} 
The quartic $\det(Q_h)=0$ is the Igusa quartic in $\PP^4=\PP(V_5^\vee)$. 
\end{prop}

\proof Recall that the generic point of the Segre cubic $C_3\subset\PP(V_5)$ is the kernel of one of the two-forms $\theta(v)$, 
and that we can get this kernel as the line generated by $\theta(v)\wedge\theta(v)\in\wedge^4V_5^\vee\simeq V_5$. 
At this generic point, the affine tangent space to the Segre cubic is therefore the hyperplane of $V_5$ generated 
by the vectors of the form $\theta(v)\wedge\theta(w)\in\wedge^4V_5^\vee\simeq V_5$. This hyperplane is defined by a linear 
form $h_v\in V_5^\vee$ that vanishes on these vectors, which exactly means that $h_v\wedge  \theta(v)\wedge\theta(w)=0$
for any $w\in V_5$. In other words, $Q_{h_v}(\theta(v),\theta(w))=0$ for all $w\in V_5$, which means that 
$\theta(v)$ belongs to the kernel of the quadratic form $Q_{h_v}$. In particular the latter is degenerate.

We have thus proved that the generic point of the projective dual variety of the Segre cubic is contained in the quartic 
hypersurface $\det(Q_h)=0$. But this projective dual is well-known to be the Igusa quartic in $\PP(V_5^\vee)$, and these
two quartics have to coincide. \qed 

\medskip\noindent 
This yields a simple determinantal representation of the Igusa quartic. Using Ozeki's representative we get
$$Q_h=\begin{pmatrix} 2h_1& -h_2 & -h_3& -h_5\\ -h_2& 2h_3 & -h_4 & 0 \\ -h_3 & -h_4& 2h_5&- h_1\\ -h_5& 0& -h_1& 2h_3  \end{pmatrix},$$
whose determinant is readily computed to be 
$$-\det(Q_h)=4h_3^4+4h_3^2(3h_1h_5-h_2h_4)-4h_3(h_1^3+h_5^3+h_1h_4^2+h_2^2h_5)+(h_1h_2-h_4h_5)^2.$$

\medskip 
One can consider inside $\PP(V_5^\vee)\times G(2,V_4)$ the locus $J_5$ of pairs $([h],U)$ such that $U$ is isotropic with respect 
to $Q_h$. Recall that $OG_Q(2,4)=\PP^1\cup\PP^1$ is the disjoint union of two smooth conics when $Q$ is non degenerate. 
When $Q$ is a quadratic form of corank one on $V_4$, the corresponding orthogonal 
Grassmannian $OG_Q(2,V_4)$ is a single conic (while if $Q$ has corank two,  $OG_Q(2,V_4)$ is the union of two planes meeting 
at one point, defined by the kernel). This means that the Stein factorization of the projection of $J_5$ to $\PP(V_5^\vee)$
is $J_5\ra Cob_4\ra \PP(V_5^\vee)$ where $Cob_4$ is the double cover of $\PP(V_5)$ branched over the Igusa quartic: that is,
the {\it Coble fourfold} \cite{cks}.

\medskip
On the other hand, denote by $\mathcal{Q}_{X_4}$ the pull-back to $X_4$ of the rank two quotient bundle on $G(3,V_5)$. 
The $\PP^1$-bundle $\PP(\mathcal{Q}_{X_4})$ over $X_4$ has a natural map to $G(2,V_4)\times \PP(V_5^\vee)$ and we claim 
that its image is precisely $J_5$. Indeed, a generic element $(U,V\subset H)$ of $G(2,V_4)\times Fl(3,4,V_5)$ belongs to 
$\PP(\mathcal{Q}_{X_4}^\vee)$ when $V$, hence $H$, contains the kernels of all the two-forms $\theta(v)$,  $v\in U$. 
But if $h$ is a linear form defining $H$, the condition that $H$ contains the kernel of $\theta(v)$ exactly means that 
$h\wedge \theta(v)\wedge \theta(v)=0$, hence our claim. 
Moreover the projection map $\PP(\mathcal{Q}_{X_4}^\vee)\ra J_5$ is birational, being clearly bijective 
outside the five special planes in $G(2,V_4)$. So we get a diagram

$$\xymatrix{& \PP(\mathcal{Q}_{X_4}^\vee)\ar[dl]_{bir}\ar[dr]^{\PP^1} & \\
J_5\ar[dr]^{conic} & & X_4\\ & Cob_4&  }$$

\medskip\noindent
where the south-west arrow is a  conic bundle, at least generically. But the picture does not seem to recover
the small resolutions of $Cob_4$ described in \cite{cks}.


\section{Local rigidity}

Since $V_4^\vee\otimes\wedge^2V_5^\vee$ is prehomogeneous, we expect that $X_4$ has strong rigidity properties. 
What we can prove is the following statement.

\begin{prop} 
$X_4$ is locally rigid and has finite automorphism group.
\end{prop} 

\proof Local rigidity is equivalent to the vanishing of  $H^1(TX_4)$. 
In order to check this, as usual we rely on the normal exact sequence, which yields an 
exact sequence of cohomology groups 
$$H^0(TG_{|X_4})\lra H^0(E_{|X_4})\lra H^1(TX_4)\lra H^1(TG_{|X_4}).$$
So local rigidity will follow from the following statements, to be proved separately:
\begin{enumerate}
\item $H^1(TG_{|X_4})=0$;
\item $H^0(E_{|X_4})=V_4^\vee\otimes\wedge^2V_5^\vee/\langle\theta\rangle$;
\item $H^0(TG_{|X_4})\lra H^0(E_{|X_4})$ is surjective.
\end{enumerate}
The third statement follows from the fact that $\PP(V_4^\vee\otimes\wedge^2V_5^\vee)$ is prehomogeneous under 
$PGL(V_4)\times PGL(V_5)$, more precisely from the fact that the orbit of $[\theta]$ is open, since this implies that 
the image of the natural differential $\fsl_4\times\fsl_5=H^0(TG)\lra V_4^\vee\otimes\wedge^2V_5^\vee/\langle\theta\rangle$
sending $X$ to $X(\theta)$ mod $\theta$ is surjective; since this morphism can also be defined by restricting first to $X_4$ and then
composing with the morphism we are interested in, the latter must also be surjective. 

In order to prove the second statement we twist by $E$ the Koszul complex resolving the structure sheaf of $X_4$. 
By standard cohomological arguments, it is enough to check that $H^0(G,End_0(E))=0$ and 
$H^i(G,E\otimes\wedge^{i+1}E^\vee)=0$ 
for any $i>0$. For the first claim, observe that $$End_0(E)=End_0(U)\oplus End_0(V)\oplus End_0(U)\otimes End_0(V)$$
is in fact acyclic. For the second claim, check that $E\otimes\wedge^{i+1}E^\vee$ is also acyclic for any $i>0$. 
Similarly, in order to prove the first statement we need to check that $H^{i+1}(G,TG\otimes\wedge^{i}E^\vee)=0$
for any $i\ge 0$, which is again a straightforward application of Bott's theorem.  

By the same type of arguments (or using a computer to check that $\chi(TX_4)=0$), we deduce that $H^0(TX_4)=0$, 
which implies that the automorphism group is discrete, hence finite since $X_4$ is Fano. 
\qed

\medskip The question remains open, whether $X_4$ is also globally rigid, which would be remarkable 
for a Fano fourfold with such a big Picard number. The first thing to check is whether $X_4$ remains smooth when 
we degenerate $\theta$ to the codimension one orbit. If yes, we would get a similar situation to the case of 
codimension two linear sections of the spinor tenfold (which has Picard number one, though). 

Another question one may ask is whether the quotient bundle restricted to 
$X_4$ is rigid. In other words, is the morphism to  $G(3,V_5)$ uniquely 
defined?

\section{Higher dimensions} 

Let us briefly describe the higher dimensional models.

\begin{prop} 
$X_6$ is a rational Fano sixfold of index one and Picard rank two.

The projection of $X_6$ to $G(3,V_5)$ is birational, with non trivial fibers isomorphic
to $\PP^1$ over the smooth locus of $C_4$, and to $\PP^2$ over its ten singular points.

The projection to $\PP(V_4)$ is a $\QQ^3$-bundle outside $P_0$, with five four-dimensional
fibers over $P_0$.
\end{prop}

From this description and that of $X_4$, we deduce that in the Grothendieck ring of varieties one has the relation
$[X_6]+L^3[Y_0]=[G(3,V_5)]+L[X_4]$. This yields the Poincaré polynomial of $X_6$, 
$$P_{X_6}(t)=1+2t+8t^2+9t^3+8t^4+2t^5+t^6.$$

\begin{prop} 
$X_8$ is a Fano eightfold of pseudo-index three, while $X'_8$ is Fano eigthfold of index three. 
\end{prop} 


\begin{prop} 
The projections of $X_8, X'_8$ to $G(2,V_5)$ are dual $\PP^2$-fibrations over the complement of a del Pezzo surface of degree five, the exceptional fibers 
being isomorphic to $\PP(V_4)$ and $G(2,V_4)$ respectively. 
\end{prop} 

We can readily deduce that $X_8$ and $X'_8$ have pure cohomology, with Poincaré polynomials
$$\begin{array}{lll}
P_{X_8}(t) & = & 1+2t+4t^2+6t^3+11t^4+6t^5+4t^6+2t^7+t^8, \\
P_{X'_8}(t) & = & 1+2t+5t^2+11t^3+13t^4+11t^5+5t^6+2t^7+t^8.
\end{array}$$

\medskip Of course $X_6, X_8, X'_8$  inherit the same symmetries as $X_4$. E. Fatighenti and F. Tanturri checked the 
necessary vanishing conditions to establish, as for $X_4$, that they are also infinitesimally rigid.

    \bibliographystyle{amsplain}

\bibliography{segre.bib}

\end{document}